\documentclass{amsart}
\usepackage{latexsym}
\usepackage{amssymb,amsthm,amsmath}

\theoremstyle{plain}
\newtheorem{theorem}{Theorem}

\newtheorem{corollary}[theorem]{Corollary}

\newtheorem{example}[theorem]{Example}

\theoremstyle{definition}
\newtheorem{definition}[theorem]{Definition}
\theoremstyle{remark}
\newtheorem{remark}[theorem]{Remark}

\newcounter{zd}

\newcounter{zdr}[subsection]

\newcommand{\eps}{\varepsilon}
\newcommand{\pa}{\partial}
\def\pa{\partial}
\let\cal=\mathcal
\begin{document}

\title[Delta shock wave formation]{Delta shock wave formation in the case of
triangular hyperbolic system of conservation laws}\thanks{The work of the first
author is supported be RFFI grant 05-01-00912,
 DFG Project 436 RUS 113/895/0-1.}
\author{
V.G.~Danilov}\address{ Moscow Technical University of Communication
and Informatics, Moscow,
Russia}\email{danilov@miem.edu.ru}\author{D.~Mitrovic}\address{Faculty
of Mathematics and Natural Sciences, University of Montenegro, 81000
Podgorica, Montenegro}\email{matematika@cg.yu}
\date{}

\maketitle

\begin{abstract}
We describe $\delta$ shock wave arising from continuous initial data
in the case of triangular conservation law system arising from
"generalized pressureless gas dynamics model" (e.g. \cite{huang2}).
We use the weak asymptotic method \cite{DSO,DM}.
\end{abstract}

In this paper we investigate formation of $\delta$-shock wave in the
case of triangular system of conservation laws:
\begin{align}
&u_t+(f(u))_x=0,
\label{jul767}\\
&v_t+(v g(u))_x=0,
\label{jul361}
\end{align} with continuous initial data
\begin{align}
u|_{t=0}=\hat{u}(x)&=\begin{cases}
U_1, \ \ x<a_2,\\
u_0(x), \ \ a_1\leq x \leq a_2, \\
U_0, \ \ a_1< x
\end{cases},
\label{jul769}\\
v|_{t=0}=\hat{v}(x)&=\begin{cases}
V_1, \ \ x<a_2,\\
v_0(x), \ \ a_1\leq x\leq a_2, \\
V_0, \ \ a_1< x
\end{cases}
\label{jul362}
\end{align} where $u_0$ and $v_0$ are continuous functions defined on $[a_2,a_1]$ such that $v_0$ is bounded and $u_0$ satisfies:
$$f'(u_0(x))=-Kx+b, \ \ x\in [a_2,a_1],$$ and $K$ and $b$ are constants determined from the continuity conditions:
$$f'(U_1)=-Ka_2+b, \ \ f'(U_0)=-K a_1+b$$ i.e.
$$K=\frac{f'(U_1)-f'(U_0)}{a_1-a_2}, \\ b=\frac{f'(U_1)a_1-f'(U_0)a_2}{a_1-a_2}.$$

For the functions $f$ and $g$ we assume:
\begin{align*}
&f\in C^2([U_0,U_1]), \ \ g\in C^1([U_0,U_1]),\\
&f''>0 \text{   on   } [U_0,U_1],\\
&g'-f''\geq 0 \text{   on   } [U_0,U_1],\\
&\exists \hat{U}\in (U_0,U_1) \text{   such that   }
g(\hat{U})=f'(\hat{U}).
\end{align*} As we will see in the next section, such conditions provide appearance of admissible $\delta$ shock wave.

The main result of the paper is construction of formulas which
smoothly and globally in $t\in {\bf R}$ represents approximate
solution to the problem. We have chosen initial data so that such
formulas correspond to the process of delta shock wave formation
from continuous initial data.

It is well known that if the solution to equation (\ref{jul767}) is
discontinuous function then unknown function $v$ contained in
(\ref{jul361}) can contain $\delta$ distribution. Such form of the
function $v$ is natural from the viewpoint of applications.

The basic difficulty lies in giving the sense to the solution
containing $\delta$ distribution. Namely, the second equation of the
system contains nonlinearity generally implying the problem of defining the product
of $\delta$ distribution with Heaviside function.

Of course, this problem appears only if we directly substitute
functions containing mentioned singularities into the equation. In
that case, the product of $\delta$ and Heaviside function is defined
using measure theory \cite{huang1,huang2,chi,ding}. Also, one can
use regularization of $\delta$ and Heaviside distribution and define
the product as the weak limit of the product of the approximations
\cite{keyf,MN}. Finally, recently it has been proposed \cite{DSH1}
to define solutions of $\delta$ type (i.e. solutions containing
$\delta$ distributions) trough appropriate integral equalities as
done in the case of $L^{\infty}\cap L^1$ weak solutions, thus
avoiding problem of multiplication of singularities (see also
\cite{huang3}).

According to all said above, we see that the problem of propagation
of already formed $\delta$ shocks has been explored rather
thoroughly.

But, the problem of $\delta$ shock formation is much less studied.
Even if it is studied, it has been always done for Riemann problem
and always using vanishing viscosity approach. One of the first
result of finding global smooth approximation to the problem of type
(\ref{jul767})-(\ref{jul362}) can be found in \cite{ind} for the
case $f(u)=u^2/2$ and $g(u)=u$. There, approximate solution $u_{\eps}$ to the
first equation of the system is found in explicit form by using vanishing viscosity
approximation and Hopf-Cole transformation. Then, substituting
$u_{\eps}$ in the place of $u$ in equation (\ref{jul361}), the equation becomes linear equation in $v$
and it is solved by the method of characteristics.

General situation using the same approach (vanishing viscosity but
with the vanishing term of the form $\eps t (u,v)_{xx}$) was
considered in \cite{ercole}. There, for the Riemann initial data
author proves that system (\ref{jul767}),(\ref{jul361}) admits
approximate solutions converging to $\delta$ type distribution.
Author obtains the result by using various apriori estimates
obtained from the equations with vanishing viscosity.

Here, we consider more general initial conditions (continuous
initial conditions), and also give explicit formula for approximate
solution to the problem. It is very important when generalizing
result to multidimensional situation since it provides us to describe
analytically geometric singularities \cite{dan4}.

Other possible approach in general situation is to use Oleinik-Lax
formula which gives the solution of the first equation
\cite{huang1}. Still, in this case one has to add additional
assumptions concerning behavior of the solution in the moment of
bifurcation.

Assumptions on the weak continuity of the solution are not enough
\cite{dsh3}, and additional assumptions on boundedness are not quite
justified for $\delta$ shock wave type solutions.

In this paper, we will use the weak asymptotic method which appeared
to be rather effective in lots of situations involving nonlinear
waves formation and interaction \cite{DSH,DSH1,DSO,DM}.

We give basic definitions of the weak asymptotic method.

\begin{definition}
By $O_{{\cal D}'}(\eps^\alpha)\subset {\cal D}'({\bf R})$,
$\alpha\in {\bf R}$, we denote the family of distributions depending
on $\eps\in (0,1)$ and $t\in {\bf R}^+$ such that for any test
function $\eta(x)\in C^1_0({\bf R})$, the estimate
$$
\langle O_{{\cal D}'}(\eps^\alpha),\eta(x)\rangle=O(\eps^\alpha), \
\ \eps\to 0,
$$
holds, where the estimate on the right-hand side is understood in
the usual sense and locally uniform in~$t$, i.e.,
$|O(\eps^\alpha)|\leq C_T\eps^\alpha$ for $t\in[0,T]$.
\end{definition}

\begin{definition}
The family of pairs of functions
$(u_\eps,v_{\eps})=(u_\eps(x,t),v_\eps(x,t))$, $\eps>0$, is called a
weak asymptotic solution of problem~(\ref{jul767})-(\ref{jul362}) if
\begin{align*}
&u_{\eps t}+(f(u_{\eps}))_x={\cal O}_{{\cal D}'}(\eps),\\
&v_{\eps t}+(v_{\eps} g(u_{\eps}))_x={\cal O}_{{\cal D}'}(\eps),\\
& u_\eps\bigg|_{t=0}-\hat{u}={\cal O}_{{\cal D}'}(\eps), \ \
v_{\eps}\bigg|_{t=0}-\hat{v}={\cal O}_{{\cal D}'}(\eps), \ \ \eps\to
0.
\end{align*}
\label{jul3126}
\end{definition}

According to the previous definitions, an approximation constructed
by the means of the vanishing viscosity is indeed weak asymptotic
solution to equation (\ref{jul767}). In the case of the quadratic
nonlinearity (i.e. when $f(u)=u^2$) weak asymptotic solution is
constructed in \cite{dan}, and in the case when $f$ is arbitrary
convex function in \cite{DM}. In both of the latter papers Cauchy
problems with special initial data of type (\ref{jul769}) are
considered.

Such initial data has property that the solution of corresponding
Cauchy problem can always be represented as linear combination of
the Heaviside functions (see \cite{DM} and Theorem 6 below). This,
in turn, allows us to use the techniques of the weak asymptotic
method \cite{DSH,DSO}.

We stress that in \cite{dan, DM}, in the case of
gradient catastrophe the solution is approximated by
the smooth function which describes interaction of nonlinear waves.

Still, in the case of general piecewise monotone continuous initial
data the problem can be solved by the same technique as in the case
of initial data (\ref{jul769}). We describe briefly the procedure
\cite{dm1}.

On the first step, we replace given (arbitrary piecewise monotone
continuous) initial data by the new initial data which differs from
the original one only in the neighborhood
$(x_0-\eps^{1/2},x_0+\eps^{1/2})$ of the point $x_0$ from which
emanates characteristic along which the shock wave will appear the
first (i.e. along which we will have the first gradient
catastrophe). In the interval $(x_0-\eps^{1/2},x_0+\eps^{1/2})$ the
new initial data will be of the type (\ref{jul769}) (where
$a_2=x_0-\eps^{1/2}$ and $a_1=x_0+\eps^{1/2}$). The discrepancy
caused by such replacement is of order ${\cal O}(\eps^{1/2})$ in
$L^1({\bf R})$ sense.

With new initial data we can construct the solution using formulas
from \cite{DM} or Theorem 6 below.

In this paper, we give new formulas for the construction of the weak
asymptotic solution to problem (\ref{jul767}), (\ref{jul769}) (which
are simpler then ones appearing in \cite{DM}). Then, after obtaining
the weak asymptotic solution to (\ref{jul767}), (\ref{jul769}) we
substitute it in equation (\ref{jul361}) and solve problem
(\ref{jul361}), (\ref{jul362}) using the method of characteristics.
We will explain this more closely.

The system we consider is triangular system. The first equation is
scalar conservation law with convex nonlinearity. The problem of
shock wave formation for this equation is solved in \cite{DM}. It
was done trough the concept of 'new characteristics'. In a matter of
fact, generalized characteristics (see \cite{daf}, Definition
10.2.1) are approximated by 'new characteristics'. But, unlike
generalized characteristics which are nonsmooth and intersecting
curves, the 'new characteristics' are globally defined smooth
nonintersecting curves along which solution remains constant. Using
this concept, we obtain approximating solution (in the sense of
Definition \ref{jul3126}) to equation  (\ref{jul767}). This solution
is smooth and we can substitute it into equation (\ref{jul361}) and
differentiate according to Leibnitz rule the product $v g(u)$
appearing there. Then, we can solve obtained equation using standard
method of characteristics.

According to the weak asymptotic method and the procedure previously
described, we replace problem (\ref{jul767}-\ref{jul362}) by the
family of problems
\begin{align}
u_{\eps t}+(f(u_{\eps})_x&={\cal O}_{{\cal D}'}(\eps),
\label{okt436}\\
v_{\eps t}+(v_{\eps}g(u_{\eps})_x&=0 \ \ \eps>0, \label{okt446}
\end{align}
\begin{align}
u_{\eps}|_{t=0}=\hat{u}(x)+{\cal O}_{{\cal D}'}(\eps)&=\begin{cases}
U_1, \ \ x<a_2,\\
u_0(x), \ \ a_1\leq x<a_2, \\
U_0, \ \ a_2\leq x
\end{cases}+{\cal O}_{{\cal D}'}(\eps),
\label{okt456}\\
v_{\eps}|_{t=0}=\hat{v}(x)+{\cal O}_{{\cal D}'}(\eps)&=\begin{cases}
V_1, \ \ x<a_2,\\
v_0(x), \ \ a_1\leq x<a_2, \\
V_0, \ \ a_2\leq x
\end{cases}+{\cal O}_{{\cal D}'}(\eps).
\label{okt466}
\end{align} where ${\cal O}_{{\cal D}'}$ will be precised in Theorem 6.

We expose the plan of the paper in more details.

In Section 1 we recall necessary conditions for appearance of
admissible $\delta$ shock wave for system (\ref{jul767}),
(\ref{jul361}). Then, we quote result in the framework of the
weak asymptotic method that we shall need.

In Section 2 we construct the weak asymptotic solution to problem
(\ref{jul767}), (\ref{jul769}).

In Section 3 we construct the weak asymptotic solution to problem
(\ref{jul361}), (\ref{jul362}).

Finally, in Section 4 we find weak limit of the constructed weak
asymptotic solution to problem (\ref{jul767}-\ref{jul362}).

\section{Conditions for $\delta$ shock wave appearance and some weak asymptotic formulas}

Consider system (\ref{jul767}), (\ref{jul361})  with Riemann initial
data:
\begin{align}
u|_{t=0}&=\begin{cases}
U_l, \ \ x<0,\\
U_r, \ \ x\geq 0,
\end{cases},\label{jul768}\\
v|_{t=0}&=\begin{cases}
V_l, \ \ x<0,\\
V_r, \ \ x\geq 0.
\end{cases}
\label{jul363}
\end{align} Since the aim of the paper is to describe formation of $\delta$ shock waves,
we want to determine sufficient condition on $f$ and $g$ which provides $\delta$-shock wave formation from initial data (\ref{jul768}), (\ref{jul363}).
In other words, we want to determine conditions on $f$ and $g$ such
that Riemann problem  (\ref{jul767}), (\ref{jul361}),
(\ref{jul768}), (\ref{jul363}) admits solution of the type:

\begin{align}
u(x,t)=&\begin{cases}
U_l, \ \ x<ct,\\
U_r, \ \ x\geq ct,
\end{cases}\\
v(x,t)=&\begin{cases}
V_l, \ \ x<ct,\\
V_0, \ \ x\geq ct
\end{cases}+const.\cdot t\cdot \delta(x-ct).
\label{jul367}
\end{align} The solution is understood in the sense of Definition 1 from paper \cite{DSH}.

As the admissibility conditions for $\delta$ shocks we shall use
overcompresive conditions (as in \cite{ercole, keyf, liu, tan}):
\begin{gather}
\lambda_i(U_r,V_r)\leq c \leq \lambda_i(U_l,V_l), \ \ i=1,2,
\label{jul364}
\end{gather}where $\lambda_i$, $i=1,2$, are eigenvalues of system (\ref{jul361}), i.e.
$$\lambda_1(u,v)=f'(u), \ \ \lambda_2(u,v)=g(u).$$
From (\ref{jul364}) and expressions for $\lambda_i$, $i=1,2$ we have:
\begin{gather}
f'(U_r)\leq c \leq f'(U_l)\notag\\
g(U_r)\leq c \leq g(U_l).
\label{jul365}
\end{gather}
The following conditions were used in \cite{ercole}:
$$g'>0, \ \ f''>0, \ \ f'<g.$$ Still, such conditions will not necessarily give $\delta$ shock
even if the classical solution $u$ to (\ref{jul767}), (\ref{jul769}) blows up after certain time. Since in this paper we are interested
only on the $\delta$ shock appearance phenomenon, we shall need more restrictive conditions. The conditions which we shall derive below ensures
$\delta$ shock wave appearance if the classical solution to (\ref{jul767}), (\ref{jul769}) blows up. We stress that $\delta$ shock wave can arise also
in the case of less restrictive conditions on $f$ and $g$ but in the special case of initial data.

We proceed with deriving of necessary conditions.
Initial assumption is convexity of the function $f$, i.e. $f''>0$.
We have to find conditions on $g$ such that (\ref{jul365}) is
satisfied. The following condition obviously implies (\ref{jul365}):
\begin{equation}
g(U_r)\leq f'(U_r)\leq c \leq f'(U_l)\leq g(U_l).
\label{jul366}
\end{equation} Since $f'$ is
increasing it is clear that it has to be $U_r>U_l$. If we assume
that $F=g-f'$ is increasing in the interval $[U_r,U_l]$ and that $F$
attains zero in that interval, obviously (\ref{jul366}) will be
satisfied (since $F$ changes sign on $[U_r,U_l]$). We can collect
previous considerations in the following theorem:
\begin{theorem}
\label{th3}
Assume that the functions $f,g\in C^2({\bf R})$ satisfy
\begin{enumerate}
\item $f''>0$ on ${\bf R}$

\item $g'-f''>0$ on $[U_r,U_l]$,

\item $\exists \hat{U}\in (U_r,U_l)$ such that $g(\hat{U})=f'(\hat{U})$.
\end{enumerate}
If $\hat{U}\in [U_r,U_l]$ then Riemann problem
(\ref{jul767}),(\ref{jul361}), (\ref{jul366}), (\ref{jul768}) admits
$\delta$ type solution of the form (\ref{jul367}). \label{jul368}
\end{theorem}

Next, we give very important theorem in the framework of the
weak asymptotic method (sometimes called nonlinear superposition
law):

\begin{theorem}\cite{DM} Let $\theta_{i\eps}(x)=\omega_i(x/\eps)$, $i=1,2$, where $\lim\limits_{z\to +\infty}\omega_i(z)=1$,
$\lim\limits_{z\to -\infty}\omega_i(z)=0$ and
$\frac{d\omega(z)}{dz}\in {\cal S}({\bf R})$ where ${\cal S}({\bf
R})$ is the Schwartz space of rapidly decreasing functions. For the
bounded functions $a,b,c$ depending on $(x,t)\in {\bf R}^+\times
{\bf R}$ we have
\begin{multline}
\label{an0}
f\left(a+b\theta_{1\eps}(\varphi_1-x)+c\theta_{2\eps}(\varphi_2-x)\right)=\\
f(a)+\theta_{1\eps}(\varphi_1-x) \left(
f(a+b+c)B_1+f(a+b)B_2-f(a+c)B_1-f(a)B_2\right)+\\
\theta_{2\eps}(\varphi_2-x)
\left(f(a+b+c)B_2-f(a+b)B_2+f(a+c)B_1-f(a)B_1\right)+ {\cal
O}_{{\cal D}'}(\eps),
\end{multline}where for $\rho\in {\bf R}$ we have
\begin{gather}
B_1(\rho)= \int\dot{\omega}_1(z)\omega_2(z+\rho)dz \text{  and }
B_2(\rho)= \int\dot{\omega}_2(z)\omega_1(z-\rho)dz, \label{1}
\end{gather} and
$$B_1(\rho)+B_2(\rho)=1.$$
\end{theorem}

\section{Weak asymptotic solution to (\ref{jul767}), (\ref{jul769})}

First, we determine the function $u_{\eps}$.  In the sequel we use the following notation (as usual $x\in {\bf R}$, $t\in {\bf R}^+$):

\begin{align*}
u_1&=u_1(x,t,\eps), \ \ B_i=B_i(\rho), \ \
\varphi_i=\varphi_i(t,\eps),\nonumber\\
\theta_{i\eps}&=\theta_{i\eps}(\varphi_i-x)=\omega_i(\frac{\varphi_i-x}{\eps}),
\\
\delta_{i\eps}&=-\frac{d}{dx}\theta_{i\eps}(\varphi_i-x)=-\frac{d}{dx}\omega_i(\frac{\varphi_i-x}{\eps}),
\; i=1,2,\\
\tau&=\frac{f'(U)t+a_2-f'(u_0^0)t-a_1}{\eps}=\frac{\psi_0(t)}{\eps},\\
t^*&=\frac{a_1-a_2}{f'(U_1)-f'(U_0)}, \\
x^*&=f'(U_1)t^*+a_2=f'(U_0)t^*+a_1=\frac{f'(U_0)a_1-f'(U_0)a_2}{f'(U_1)-f'(U_0)}.
\end{align*} The function $\tau$ is so-called 'fast variable'. It is equal to difference of standard characteristics
of equation (\ref{jul361}) emanating from $a_2$ and $a_1$,
respectively. When we are in the domain of existence of classical
solution to (\ref{jul767}), (\ref{jul769}) we have $\tau\to
-\infty$, while when we are in the domain where solution to
(\ref{jul767}, (\ref{jul769}) is discontinuous (i.e. in the form of
the shock wave) we have $\tau \to \infty$.

The point $(t^*,x^*)$ is the point of blow up of the classical solution to (\ref{jul767}), (\ref{jul769}).

Also notice that $w-\lim\limits_{\eps\to
0}\theta_{i\eps}=\theta(\varphi_i-x)$ and $w-\lim\limits_{\eps\to
0}\delta_{i\eps}=\delta(\varphi_i-x)$ for Heaviside function
$\theta$ and Dirac distribution $\delta$.

First, we will describe the shock wave formation process for problem
(\ref{jul767}), (\ref{jul769}). We explain the procedure we shall
use before we formulate the theorem.

It is well known problem (\ref{jul767}), (\ref{jul769}) will have
classical solution up to the moment $t=t^*$. The choice of our
initial data is such that in the moment of blow up of the classical
solution the shock wave will be formed and it will not change its
shape for any $t>t^*$. This is because all the characteristics
emanating from $[a_2,a_1]$ intersect in one point $(t^*,x^*)$.


So, for $t>t^*$ we have to pass to the weak solution concept. In
other words, in the moment $t=t^*$ we stop the time and solve
Riemann problem for equation (\ref{jul767}).

Our aim is to find global in time approximate solution to
(\ref{jul767}), (\ref{jul769}) which is at least continuous. To do
this we have to avoid intersection of characteristics.

Natural idea is to smear the discontinuity line, i.e. to take $\eps$
neighborhood of the discontinuity line and to dispose
characteristics in that neighborhood in a way that they do not
intersect and as $\eps \to 0$ all of them lump together into the
discontinuity line. Of course, this will not be the standard
characteristics for problem (\ref{jul767}), (\ref{jul769}).
Nevertheless, along them approximate solution to our problem will
remain constant. Such lines we call 'new characteristics'.

Another question that arises here is how to distribute 'new
characteristics' in the $\eps$ neighborhood of the discontinuity
line. The obvious way to accomplish this is to distribute the 'new
characteristics' uniformly in the mentioned area, i.e. in a way that
every of them is parallel to the discontinuity line.

Since all the characteristics emanating from the interval
$[a_2,a_1]$ intersect in the same point, roughly speaking, it is
enough to find the way to dispose 'new characteristics' emanating
from $a_2$ and $a_1$ so that they do not intersect.

We use Theorem \ref{DM} and 'switch' functions $B_i$, $i=1,2$,
appearing there.

Denote by $\varphi_i$, $i=1,2$, the new characteristics emanating
from the points $a_i$, $i=1,2$, respectively. They are given by the
following Cauchy problems:
\begin{equation}
\begin{split}
&\frac{d}{dt}\varphi_{1}(t,\eps)=(B_2(\rho)-B_1(\rho))f'(U_1)+cB_1(\rho),
\ \ \varphi_{1}(0,\eps)=a_1+A\eps \frac{a_1+a_2}{2},\\
&\frac{d}{dt}\varphi_{2}(t,\eps)=(B_2(\rho)-B_1(\rho))f'(U_0)+cB_1(\rho),
\ \ \varphi_{2}(0,\eps)=a_2-A\eps \frac{a_1+a_2}{2},
\end{split}
\label{trnd1}
\end{equation} for large enough constant $A$. As we shall see later, it will be necessary to extend a little bit the interval $[a_2, a_1]$.
Therefore, we have $A\eps \frac{a_1+a_2}{2}$ accompanying initial
data in (\ref{trnd1}). Also, in (\ref{trnd1}) we define:
\begin{gather}
\rho=\frac{\varphi_2(t,\eps)-\varphi_1(t,\eps)}{\eps}.
\label{no1}
\end{gather}
According to what we said above, we expect that for every $t>t^*$ it
should be (since new characteristics should be 'close' one to
another for $t>t^*$, i.e. equal as $\eps\to 0$):
$$\varphi_1(t,\eps)-\varphi_2(t,\eps)={\cal O}(\eps), \ \ t>t^*,$$ and (since
new characteristics should be 'close' to the discontinuity line):
$$\frac{d}{dt}\varphi_{i\eps}(t,\eps)-\frac{f(U_1)-f(U_0)}{U_1-U_0}={\cal
O}(\eps), \ \ t>t^*.$$ More precisely, we expect that for $t>t^*$
the expression $B_2(\rho)-B_1(\rho)$ to be close to zero thus
eliminating nonlinearity $f'$ appearing in the equation of new
characteristics (\ref{trnd1}). This means that, according to Theorem
\ref{DM} we have $B_2(\rho)+B_1(\rho)=1$, $B_1$ and $B_2$ are close
to $1/2$, and this implies that $c$ from (\ref{trnd1}) should be
close to $2\frac{f(U_1)-f(U_0)}{U_1-U_0}$ (Rankine-Hugoniot
conditions). Actually, here we use the following simple observation.

Once the shock wave is formed, it continuous to move according to
Rankine-Hugoniot conditions and it does not change its shape along
entire time axis. Therefore, the linear equation:
\begin{equation}
\frac{\pa u}{\pa t}+\frac{c}{2}\frac{\pa u}{\pa x}=0, \ \
c=2\frac{f(U_1)-f(U_0)}{U_1-U_0}, \label{jul161}
\end{equation}
and equation (\ref{jul767}) with the same initial condition:
\begin{equation*}
u|_{t=0}=\begin{cases}
U_1, \ \ x<0,\\
U_0, \ \ x\geq 0,
\end{cases}
\end{equation*} will have the same solutions. Clearly, it is much
easier to solve linear equation (\ref{jul161}) then nonlinear
equation (\ref{jul767}). Still, the question is how to pass from
nonlinear equation (\ref{jul767}) to linear equation (\ref{jul161})
in the domains where they give the same solution (in the case of our
initial data it will be after the shock wave formation). We explain
briefly how we do it.

Define 'new characteristics'  as the solutions to the following
Cauchy problem:
\begin{equation}
\begin{split}
\dot x&=(B_2-B_1)f'(u_1)+cB_1,\\
\dot u_1&=0, \ \ u_1(0)=u_0(x_0),  \ \ x(0)=x_0+\eps A
(x_0-\frac{a_1+a_2}{2}), \ \ x_0 \in [a_2,a_1].
\end{split}
\label{avg45}
\end{equation} Thus, $\varphi_i(t,\eps)=x(a_i,t,\eps)$, $i=1,2$,
where $x$ is the solution to (\ref{avg45}). We will show later that
it is possible to choose the constant $A$ so that for $x_0\in
[a_2,a_1]$ every $t>0$ we have
$$
\frac{\pa x}{\pa x_0}>0.
$$ This means that 'new characteristics' indeed do not intersect
which in turn means that there exists the solution $x_0$ of the
implicit equation:
\begin{gather}
\label{trnd103}
x(x_0,t,\eps)=x.
\end{gather} Bearing in mind that $B_1\sim B_2$ after the
interaction we see that, using the new characteristics, we have
smoothly passed from the characteristics of equation (\ref{jul767})
to the characteristics of equation (\ref{avg45}), i.e. form equation
(\ref{jul767}) to equation (\ref{jul161}).


We formalize the previous considerations in Theorem \ref{t**}. The
theorem is analogue to the main result from \cite{DM}. The problem
which we consider here, i.e. problem (\ref{jul767}), (\ref{jul769})
can be solved in more elegant manner (see Theorem \ref{tjul861}
below). Still, approach used in Theorem \ref{t**} can be used on the
case of arbitrary piecewise monotone initial data. Also, Theorem
\ref{t**} represents motivation for Theorem \ref{tjul861}.

\begin{theorem}
\label{t**}
The weak asymptotic solution of problem (\ref{jul767}), (\ref{jul769}) has the form:
\begin{multline}
u_{\eps}(x,t)=U_0+\left(u_1(x,t,\eps)-U_0\right)\omega_{1}(\frac{\varphi_1(t,\eps)-x}{\eps})+\\
\left(U_1-u_1(x,t,\eps)\right)\omega_{2}(\frac{\varphi_2(t,\eps)-x}{\eps}),
\label{avg725}
\end{multline}where $\omega_i$ satisfies the conditions from Theorem 5 and:

$$\omega_1(z)=1 \ \  for \ \ z>0 \ \  and \ \ \omega_2(z)=0 \ \  for \ \ z<0.$$

The functions $\varphi_{i}(t,\eps)$, $t\in {\bf R}^+$, $i=1,2$, are
given by (\ref{trnd1}) and the function $\rho$ is given by
(\ref{no1}).


The function $u_1(x,t,\eps)$ is given by
$$u_1(x,t,\eps)=u_0(x_0(x,t,\eps))$$ where $x_0$ is the inverse
function to the function $x=x(x_0,t,\eps)$, $t>0$, $\eps>0$, of 'new
characteristics' defined trough Cauchy problem  (\ref{avg45}).
\end{theorem}
{\bf Proof:} We substitute anzatz (\ref{avg725}) into
(\ref{okt436}):
\begin{multline*}
\Big(U_0+\left(u_1(x,t,\eps)-U_0\right)\omega_{1}(\frac{\varphi_1(t,\eps)-x}{\eps})+
\left(U_1-u_1(x,t,\eps)\right)\omega_{2}(\frac{\varphi_2(t,\eps)-x}{\eps})
\Big)_t+\\\Big(
f(U_0+\left(u_1(x,t,\eps)-U_0\right)\omega_{1}(\frac{\varphi_1(t,\eps)-x}{\eps})+
\left(U_1-u_1(x,t,\eps)\right)\omega_{2}(\frac{\varphi_2(t,\eps)-x}{\eps}))\Big)_x\\={\cal
O}_{{\cal D}'}(\eps)
\end{multline*}
After using (\ref{an0}) and Leibnitz rule, and collecting terms
multiplying
$\theta_{1\eps}=\omega_{1}(\frac{\varphi_1(t,\eps)-x}{\eps})$,
$\theta_{2\eps}=\omega_{2}(\frac{\varphi_2(t,\eps)-x}{\eps})$,
$\delta_{1\eps}=\left(\omega_{1}(\frac{\varphi_1(t,\eps)-x}{\eps})\right)'_x$
and
$\delta_{2\eps}=\left(\omega_{2}(\frac{\varphi_2(t,\eps)-x}{\eps})\right)'_x$
we have:
\begin{align*}
&\left[\frac{\pa u_1}{\pa t}+B_2(\rho) f'(u_1)\frac{\pa u_1}{\pa x}+
B_1(\rho) f'(U_1+U_0-u_1)\frac{\pa u_1}{\pa x} \right]\theta_{1\eps}+\\
&\left[-\frac{\pa u_1}{\pa t}-B_2(\rho) f'(u_1)\frac{\pa u_1}{\pa
x}-B_1(\rho)
f'(U_1+U_0-u_1)\frac{\pa u_1}{\pa x}\right]\theta_{2\eps}+\\
&\left(
(u_1-U_0)\varphi_{1t}-B_2(\rho)\left(f(u_1)-f(U_0)\right)-B_1(\rho)\left(f(U_1)-f(U_1+U_0-u_1)
\right) \right) \delta_{1\eps}+\\
&\left(
(U_1-u_1)\varphi_{2t}-B_2(\rho)\left(f(U_1)-f(u_1)\right)-B_1(\rho)\left(f(U_1+U_0-u_1)-f(U_0)
\right) \right)\delta_{2\eps}={\cal O}_{{\cal D}'}(\eps).
\end{align*} In the sequel we write only $B_i$ instead of
$B_i(\rho)$, $i=1,2$.

We rearrange this expression using the following simple formula
$C\theta_{1\eps}+D\theta_{2\eps}=(C+D)\theta_{2\eps}+C(\theta_{1\eps}-\theta_{2\eps})$:
\begin{equation*}
\begin{split}
&\left(\frac{\pa u_1}{\pa t}+[(B_2-B_1) f'(u_1)]\frac{\pa u_1}{\pa
x}\right)\left(
\theta_{1\eps}-\theta_{2\eps} \right)+\\
&B_1 [\frac{d}{d x}\left(f(U+u_0^0-u_1)+f(u_1)\right)]\left(
\theta_{1\eps}-\theta_{2\eps} \right)+\\
&\left( (u_1-U_0)\varphi_{1t}-B_2\left(f(u_1)-f(U_0)\right)-B_1\left(f(U_1)-f(U_1+U_0-u_1) \right) \right)\delta_{1\eps}+\\
&\left((U_1-u_1)\varphi_{2t}-B_2\left(f(U_1)-f(u_1)\right)-B_1\left(f(U_1+U_0-u_1)-f(U_0)
\right)  \right)\delta_{2\eps}={\cal O}_{{\cal D}'}(\eps).
\end{split}
\end{equation*}

For an unknown constant $c$ we add and subtract the term
$cB_1\frac{\pa u_1}{\pa x}$ in the coefficient multiplying $\left(
\theta_{1\eps}-\theta_{2\eps} \right)$ and then we rewrite the
last expression in the following form:
\begin{equation}
\label{avg635}
\begin{split}
&\left(\frac{\pa u_1}{\pa t}+[(B_2-B_1) f'(u_1)+cB_1]\frac{\pa
u_1}{\pa x}\right)\left(
\theta_{1\eps}-\theta_{2\eps} \right)+\\
&B_1 [\frac{d}{d x}\left(f(U+u_0^0-u_1)+f(u_1)-c u_1
\right)]\left(
\theta_{1\eps}-\theta_{2\eps} \right)+\\
&\left( (u_1-U_0)\varphi_{1t}-B_2\left(f(u_1)-f(U_0)\right)-B_1\left(f(U_1)-f(U_1+U_0-u_1) \right) \right)\delta_{1\eps}+\\
&\left((U_1-u_1)\varphi_{2t}-B_2\left(f(U_1)-f(u_1)\right)-B_1\left(f(U_1+U_0-u_1)-f(U_0)
\right)  \right)\delta_{2\eps}={\cal O}_{{\cal
D}'}(\eps).
\end{split}
\end{equation}  We put
\begin{equation*}
\frac{\pa u_1}{\pa t}+\left[(B_2-B_1) f'(u_1)+cB_1 \right]\frac{\pa
u_1}{\pa x}=0,\ \ u_1(x,0,\eps)=u_0(x), \ \ x\in [a_2,a_1].
\end{equation*}
The system of characteristics for this problem reads:
\begin{equation}
\begin{split}
\dot x&=(B_2-B_1)f'(u_1)+cB_1,\\
\dot u_1&=0, \ \ u_1(0)=u_0(x_0),  \ \ x(0)=x_0\in [a_2,a_1].
\end{split}
\label{no10oct}
\end{equation} The aim is to prove that characteristics defined by the previous system do not intersect. It appears that it is
much easier to accomplish this if we perturb initial data for $x$ in
the previous system for a parameter of order $\eps$. More precisely,
instead of (\ref{no10oct}) we shall consider system (\ref{avg45})
(the same is done in \cite{DM}).

It is clear that such perturbation changes the solution of
(\ref{avg45}) for ${\cal O}_{{\cal D}'}(\eps)$ since initial
condition in (\ref{avg45}) is continuous.

We pass to the proof that the characteristics given by (\ref{avg45})
do not intersect. From the second equation in (\ref{avg45}) it
follows $u_1\equiv u_0(x_0)$. We substitute this into the first
equation of (\ref{avg45}) and use $f'(u_0(x_0))=-Kx_0+b$, $x_0\in [a_2,a_1]$. We have:
\begin{equation}
\dot x=(B_2-B_1)(-Kx_0+b)+cB_1, \ \ x(0)=x_0+\eps A
\left(x_0-\frac{a_1+a_2}{2}\right). \label{avg615}
\end{equation}

Out of the segment $[a_2-\eps A\frac{a_1+a_2}{2},a_1+\eps
A\frac{a_1+a_2}{2}]$ initial function is constant and we define the
solution $u_1$ of problem (\ref{avg45}) to be equal to $U_1$ on the
left-hand side of the characteristic emanating from
$a_2-A\frac{a_1+a_2}{2}$ and to be equal to $U_0$ on the right-hand
side of the characteristic emanating from $a_1+A\frac{a_1+a_2}{2}$

For the functions $\varphi_1$ and $\varphi_2$ as the characteristics
emanating from $a_1+A\frac{a_1+a_2}{2}$ and $a_2-A\frac{a_1+a_2}{2}$
respectively, we have (compare to (\ref{trnd1}))
\begin{gather}
\varphi_{1t}=(B_2-B_1)(-Ka_1+b)+cB_1,
\label{avg25}\\
\varphi_{2t}=(B_2-B_1)(-Ka_2+b)+cB_1. \label{avg35}
\end{gather}

Now, we show how to effectively determine $\rho$ given by
(\ref{no1}). We apply standard procedure (see \cite{dan,DSH,DM})).
Subtracting (\ref{avg25}) from (\ref{avg35}) we get:
\begin{equation*}
(\varphi_2-\varphi_1)_t=\eps\left(\frac{\varphi_2-\varphi_1}{\eps}\right)_t=
\eps\rho_t=(B_2-B_1)\psi_0(t).
\end{equation*}Then, passing from the
"slow" variable $t$ to the "fast" variable $\tau$ we obtain (we also use
$B_2+B_1=1$):
\begin{gather}
\rho_{\tau}=1-2B_1(\rho), \ \ \frac{\rho}{\tau}\Big|_{\tau\to
-\infty}=1. \label{3}
\end{gather}

We explain the condition $\lim\limits_{\tau\to
-\infty}\frac{\rho}{\tau}=1$. We have from (\ref{avg25}) and
(\ref{avg35})
$$
\frac{\rho}{\tau}=\frac{\int_0^t2(U-u_0^0)(B_2-B_1)dt'+a_2-a_1}{2(u-u_0^0)t+a_2-a_1}.
$$Putting $t=0$ in the previous relation we see that
\begin{equation}
\frac{\rho}{\tau}\Big|_{t=0}=1. \label{no12oct}
\end{equation} When we let $\eps\to
0$ when $t=0$ we have $\tau\to -\infty$. Therefore, from
(\ref{no12oct}) it follows
$$
\frac{\rho}{\tau}\Big|_{\tau\to-\infty}=1.
$$ This relation practically means that new characteristics
emanating from $a_i$, $i=1,2$, coincides at least in the initial
moment with standard characteristics up to some small parameter
$\eps$. Still, since $\tau\to -\infty$ for every $t<t^*$ (which
means $B_1\to 0$) we see from (\ref{avg25}) and (\ref{avg35}) that
new characteristics coincides with standard ones for every $t<t^*$
up to some small parameter $\eps$.

Next, we analyze (\ref{3}). From the standard theory of ODE we see
that $\rho\to \rho_0$ as $\tau \to +\infty$ where $\rho_0$ is
constant such that $B_1(\rho_0)=B_2(\rho_0)=1/2$. That means that
after the interaction, i.e. for $t>t^*$, we have
$$\rho=\frac{\varphi_1-\varphi_2}{\eps}=\rho_0+{\cal O}(\eps)\implies \varphi_1=\varphi_2+{\cal O}(\eps), \ \ \eps\to
0,$$or, after letting $\eps\to 0$, for $t>t^*$ we have shock wave
concentrated at (see text in front of theorem for notations):
\begin{gather}
\label{trnd*}
\varphi(t)=\lim\limits_{\eps\to 0}\varphi_i(t,\eps)=\frac{c}{2}(t-t^*)+x^*.
\end{gather}

Now, we can prove global solvability of Cauchy problem
(\ref{avg45}).

Problem (\ref{avg45}) is globally solvable if characteristics
emanating from the interval
$[a_2-A\eps\frac{a_1+a_2}{2},a_1+A\eps\frac{a_1+a_2}{2}]$ do not
intersect. To prove that we will use the inverse function theorem.
We will prove that for every $t$ we have $\frac{\pa x}{\pa x_0}>0$
which means that for every $x=x(x_0,t,\eps)$, $x_0\in [a_2,a_1]$, we have
unique $x_0=x_0(x,t,\eps)$ and we can write $u_1(x(x_0,t,\eps),t)=u_0(x_0(x,t,\eps))$.

Differentiating (\ref{avg615}) in $x_0$ and integrating from $0$ to
$t$ we obtain (we remind $B_2+B_1=1$):
\begin{equation}
\frac{\pa x}{\pa x_0}=1+\eps A -K\int_0^t(B_2-B_1)dt'=1+\eps
-K\int_0^t(1-2B_1)dt'. \label{avg625}
\end{equation}

For $t\in [0,t^*]$ we have (notice that $1-Kt^*=0$):
\begin{align*}
\frac{\pa x}{\pa x_0}&=1+\eps A-K\int_0^t dt+ K\int_0^t2B_1 dt\geq \\
&1+\eps A-K\int_0^{t*} dt+ K\int_0^{t}2B_1 dt=\eps A+K\int_0^{t}2
B_1 dt>0.
\end{align*}So, everything is correct for $t\leq t^*$.

To see what is happening for $t>t^*$, initially we estimate
$1-2B_1(\rho)$ when $\tau\to \infty$. From equation (\ref{3}) we
have (we use Taylor expansion of $B_1$ around the point
$\rho=\rho_0$):
\begin{gather*}
\rho_{\tau}=1-2B_1(\rho)=-2(\rho-\rho_0)B_1'(\tilde{\rho}),
\end{gather*}for some $\tilde{\rho}$ belonging to the interval with
ends in $\rho$ and $\rho_0$. From here we see:
$$\rho-\rho_0=(\rho(\tau_0)-\rho_0){\rm exp}(\int_{\tau_0}^{\tau}-
2B_1'(\tilde{\rho})d\tau')=(\rho(\tau_0)-\rho_0){\rm
exp}((\tau_0-\tau)2B_1'(\tilde{\rho}_1))$$ for some fixed $\rho_0\in
{\bf R}$ and $\tilde{\rho}_1\in (\rho(\tau_0),\rho(\tau))\subset
[\rho(\tau_0),\rho_0]$. We remind that  $B_1'(\tilde{\rho}_1))\geq c
>0$, for some constant $c$, since $B_1$ is increasing function and
$\tilde{\rho}_1$ belongs to the compact interval
$[\rho(\tau_0),\rho_0] $. Letting $\tau\to \infty$ we conclude that
for any $N\in {\bf N}$
$$\rho-\rho_0={\cal O}(1/\tau^N), \ \ \tau\to \infty.$$ From here we have $\rho_{\tau}={\cal O}(1/\tau^N), \ \ \tau\to \infty$, since:
\begin{equation*}
\lim\limits_{\tau\to
\infty}\frac{\rho_{\tau}}{\rho-\rho_0}=\lim\limits_{\tau\to
\infty}\frac{1-2B_1(\rho)}{\rho-\rho_0}= \lim\limits_{\tau\to
\infty}-2B_1'(\rho)=-2B'_1(\rho_0)={\rm const.}<0
\end{equation*} This, in turn, means that for every $N\in {\bf N}$ and $t>t^*$ we have
\begin{gather} 1-2B_1(\rho)=\rho_{\tau}={\cal O}(\tau^{-N})={\cal O}(\eps^N), \ \ \eps\to
\label{no2} \infty, \end{gather} since for fixed $t>t^*$ we have
$\tau=\frac{\psi_0(t)}{\eps}\to \infty$ as $\eps\to 0$.

Now we can prove that $\frac{\pa x}{\pa x_0}$ for $t>t^*$. We have
\begin{align}
\label{no60oct} &\frac{\pa x}{\pa x_0}=1+\eps
A-2K\int_0^t(1-2B_1)dt'=\\&1+\eps A-2K\int_0^{t^*}(1-2B_1)dt'-
2K\int_{t^*}^{t}(1-2B_1)dt'=\nonumber\\
&\eps A+4\int_0^{t^*}B_1dt'-2K\int_{t^*}^{t}(1-2B_1)dt'>\eps
A-2K\int_{t^*}^{t}(1-2B_1)dt'. \nonumber
\end{align}Recall that
$$B_1=B_1(\rho(\tau))=B_1(\rho(\frac{\psi_0(t)}{\eps})).$$
Consider the last term in expression (\ref{no60oct}):
\begin{gather*}
2K\int_{t^*}^{t}(1-2B_1)dt'=2K\int_{t^*}^{t}(1-2B_1(\rho(\frac{\psi_0(t')}{\eps})))dt'=\\
\left(
\begin{array}{cc} \frac{\psi_0(t')}{\eps}= z \implies&
(u-u_0^0)dt'=\eps dz;\\
t^*<t'<t\implies& 0<z<\frac{\psi_0(t)}{\eps}
\end{array}
\right)=\\2K\eps\int_0^{\frac{\psi_0(t)}{\eps}}(1-2B_1(\rho(z)))dz<\eps
2KC,
\end{gather*}where $$C=\int_0^{\infty}(1-2B_1(\rho(z)))dz<\infty,$$
since from (\ref{no2}) we know $1-2B_1(\rho(z))={\cal O}(z^{-N})$,
$z\to \infty$ and $N\in {\bf N}$ arbitrary.

Therefore, for $A$ large enough (more precisely for $A>C$)we have $\frac{\pa x}{\pa x_0}>0$
what we wanted to prove.

Next step is to obtain the constant $c$. We multiply (\ref{avg635})
by $\eta\in C^1_0({\bf R})$, integrate over ${\bf R}$ with respect
to $x$ and use (\ref{avg45}) (so, we remove the first term in
(\ref{avg635})):
\begin{gather*}
\int B_1 [\frac{d}{d x}\left(f(U_1+U_0-u_1)+f(u_1)-c u_1
\right)]\left(
\theta_{1\eps}-\theta_{2\eps} \right)\eta(x)dx+\\
\left( (u_1-U_0)\varphi_{1t}\!-\!B_2\left(f(u_1)-f(U_0)\right)\!-\!B_1\left(f(U_1)-f(U_1+U_0-u_1) \right) \right)\delta_{1\eps}+\nonumber\\
\left((U_1-u_1)\varphi_{2t}\!+\!B_2\left(f(u_1)-f(U_1)\right)\!+\!B_1\left(f(U_0)-f(U_1+U_0-u_1)
\right)  \right)\delta_{2\eps}\!=\!{\cal O}(\eps).
\end{gather*}  We apply partial integration on the first integral in the previous expression to obtain:
\begin{gather}
\label{avg55}
\int B_1 [f(U_1+U_0-u_1)+f(u_1)-c u_1]\left(
\theta_{1\eps}-\theta_{2\eps} \right)\eta'(x)dx+\\
\int \left( (u_1-U_0)\varphi_{1t}-B_2 \left(f(u_1)-f(U_0) \right)+
B_1\left(f(u_1)+f(U_0)-cu_1\right)\right) \eta(x)\delta_{1\eps}
dx+\notag
\\
\int \left(
(U_1-u_1)\varphi_{2t}\!+\!B_2\left(f(u_1)-f(U_1)\right)\!-\!B_1\left(f(u_1)+f(U_1)-c
u_1 \right)\right) \eta(x)\delta_{2\eps}dx \!=\!{\cal O}(\eps).
\notag
\end{gather} Now we use the choice of regularizations $\omega_i$, $i=1,2$. Since,
\begin{gather*}
\delta_{1\eps}=\frac{1}{\eps}\omega'_1(\frac{\varphi_1-x}{\eps})=0 \text{  for  } x<\varphi_1\\
\delta_{2\eps}=\frac{1}{\eps}\omega'_2(\frac{\varphi_2-x}{\eps})=0 \text{  for  } x>\varphi_2
\end{gather*} we have from (\ref{avg55}):
\begin{align*}
&\eps\rho B_1 \int [f(U_1+U_0-u_1)+f(u_1)-c u_1]\frac{
\theta_{1\eps}-\theta_{2\eps}}{\varphi_2-\varphi_1}\eta'(x)dx+\\
&\int_{\varphi_1}^{\infty} \left( (u_1-U_0)\varphi_{1t}-B_2
\left(f(u_1)-f(U_0) \right)+
B_1\left(f(u_1)+f(U_0)-cu_1\right)\right)
\eta(x)\delta_{1\eps} dx+\\
&\int_{-\infty}^{\varphi_2} \left(
(U_1-u_1)\varphi_{2t}\!+\!B_2\left(f(u_1)-f(U_1)\right)\!-\!B_1\left(f(u_1)+f(U_1)-c
u_1 \right)\right) \eta(x)\delta_{2\eps}dx=
\end{align*}(now, we use $u_1\equiv U_1$ for $x>\varphi_1$ and $u_1\equiv
U_0$ for $x<\varphi_2$)
\begin{align}
\label{okt476}
&\eps\rho B_1 \int [f(U_1+U_0-u_1)+f(u_1)-c
u_1]\frac{
\theta_{1\eps}-\theta_{2\eps}}{\varphi_2-\varphi_1}\eta'(x)dx+\\&
\int_{\varphi_1}^{\infty}B_1\left(2f(U_0)-c U_0 \right) \eta(x)\delta_{1\eps} dx-\nonumber\\
&\int_{-\infty}^{\varphi_2}B_1\left(-2f(U_1)+c U_1
\right)\delta_{2\eps} \eta(x)dx={\cal O}(\eps). \nonumber
\end{align}
To continue, notice that we have $|\rho B_1|<\infty$ for every
$\tau\in {\bf R}$. Namely,
\begin{equation}
\begin{split}
&|\rho B_1(\rho)|\to 0 \ \ {\rm as} \ \ \tau\to -\infty \text{ since
in that case  } B_1(\rho(\tau))\sim
B_1(\tau)\sim\frac{1}{\tau^N}\sim\frac{1}{\rho^N},\\
&|\rho B_1(\rho)|\to \rho_0B_1(\rho_0) \ \ {\rm as} \ \ \tau\to
\infty \text{ since in that case  } \rho\to \rho_0.
\end{split}
\label{trnd101}
\end{equation} This fact, together with the fact that
$\delta_{i\eps}\rightharpoonup \delta(\varphi_i-x)$, reduces
expression (\ref{okt476}) to:
\begin{align}
\label{avg655}
&B_1\left(2f(U_0)-c U_0 \right) \eta(\varphi_1)-B_1\left(-2f(U_1)+c U_1 \right)\eta(\varphi_2)={\cal O}(\eps).
\end{align}

Rewrite this expression in the following manner:
\begin{align*}
&B_1\left(2(f(U_0)-f(U_1))-c (U_0- U_1)\right)\eta(\varphi_1)+B_1 \left(-2f(U_1)+c U_1  \right)\left(\eta(\varphi_2)-\eta(\varphi_1) \right)=\\
&B_1\left(2(f(U_0)-f(U_1))-c (U_0-
U_1)\right)\eta(\varphi_1)+\\
&\eps\rho B_1(\rho)\left(-2f(U_1)+c U_1
\right)\frac{\eta(\varphi_2)-\eta(\varphi_1)}{\varphi_2-\varphi_1}=^{(\ref{trnd101})}\\
&B_1\left(2(f(U_0)-f(U_1))-c (U_0- U_1)\right)\eta(\varphi_1)={\cal
O}(\eps).
\end{align*}
From here, we see that the last relation is satisfied for
\begin{equation}
c=2\frac{f(U_1)-f(U_0)}{U_1-U_0}.
\label{2}
\end{equation}The theorem is proved.
$\Box$

\begin{remark}
In the case such as our, when $U_0$ and $U_1$ are constants, is possible to replace formula
(\ref{avg725})  by
$$u_{\eps}(x,t)=\hat{u}(x_0(x,t,\eps)),$$
where, as before, the function $x_0$ is the solution to implicit
equation (\ref{trnd103}) and $\hat{u}$ are initial data
(\ref{jul769}).

The proof of this fact obviously follows after comparing
the trajectories.
We give precise formulation in the next theorem. We leave it without proof
since it is completely analogical to the proof of the
previous theorem.
\end{remark}

The difference between the previous and the next theorem is in the form of characteristics along which we solve
our problem.

In the previous theorem, for fixed $\eps$, the weak asymptotic
solution $u_{\eps}$ to (\ref{jul767}), (\ref{jul769}) was generator
of continuous semigroup of transformations (since characteristics
intersect along $x=\varphi_i$) and in the following theorem the weak
asymptotic solution $u_{\eps}$ to (\ref{jul767}), (\ref{jul769})
forms continuous group of transformation since appropriate characteristics
do not intersect.
Still, approach from the next
theorem can be used only in the case of special initial data.


\begin{theorem}
\label{tjul861}
The weak asymptotic solution $u_{\eps}$, $\eps>0$, to Cauchy problem
\begin{equation}
u_t+(f(u))_x=0, \ \ u|_{t=0}=\hat{u}(x),
\label{jul3136}
\end{equation}
is given by
\begin{gather}
u_{\eps}(x,t)=\hat{u}(x_0(x,t,\eps)),
\label{avg2956}
\end{gather} where $x_0$ is inverse function to the function $x=x(x_0,t,\eps)$, $t>0$, $\eps>0$, of 'new characteristics' defined trough the Cauchy problem:
\begin{equation}
\begin{split}
\dot{x}&=f'(u_{\eps})(B_2(\rho)-B_1(\rho))+cB_1(\rho), \ \ x(0)=x_0+\eps A \left(x_0-\frac{a_1+a_2}{\eps} \right),\\
\dot{u_{\eps}}&=0, \ \ u_{\eps}(0)=\hat{u}(x_0), \ \ x_0\in {\bf R}.
\end{split}
\label{jul661}
\end{equation} where $A$ is large enough, the functions $B_1$ and $B_2$ are defined in Theorem \ref{tjul461}, constant $c$ is given
in (\ref{jul161}) and $\rho=\rho(\psi_0(t)/\eps)$ is the solution of
Cauchy problem (\ref{no1}).

\end{theorem}

The following corollary is obvious. It claims that the weak
asymptotic solution defined in arbitrary of the previous theorems
tends to the shock wave with the states $U_1$ on the left and $U_0$
on the right (see (\ref{trnd*})):
\begin{corollary}With the notations from the previous theorems,
for $t>t^*$ the weak asymptotic solution $u_{\eps}$ to problem
(\ref{jul767}), (\ref{jul769}) we have for every fixed $t>0$:
$$
u_{\eps}(x,t)\rightharpoonup \begin{cases} U_1, \ \
x<\frac{c}{2}(t-t^*)+x^*,\\
U_0, \ \ x>\frac{c}{2}(t-t^*)+x^*,
\end{cases}
$$where $\rightharpoonup$ means convergence in
the weak sense with respect to the real variable.
\end{corollary}

\section{The weak asymptotic solution to (\ref{jul361}), (\ref{jul362})}

At the beginning of the section, we explain some general moments.

The plan is to substitute smooth function $u_{\eps}$ given by
(\ref{avg725}) into (\ref{jul361}). Thus, we obtain equation
(\ref{okt446}). Augmented by initial data (\ref{jul362}), this
linear partial differential equation of the first order has global
differentiable solution.

But, as $\eps\to 0$ we can have discontinuities in
$v-\lim\limits_{\eps\to 0}v_{\eps}$ not only on the line on which
the shock wave of standard admissible weak solution $u$ of
(\ref{jul767}), (\ref{jul769}) is supported, and which appears for
$t>t^*$. Also, for $t<t^*$ discontinuities can appear along standard
characteristics for problem (\ref{jul767}), (\ref{jul769}) emanating
from the points $a_2$ and $a_1$. Discontinuity arises due to
non-smoothness of initial data (\ref{jul769}) in the points $a_2$
and $a_1$ (see Example \ref{sep2016}).

The situation is different if instead of (\ref{avg725}) we put in
(\ref{jul361}) in the place of of $u$ the function $u_{\eps}$ given
by (\ref{avg2956}) in Theorem 7. The function $u_{\eps}$ is not
smooth since we do not use regularizations of Heaviside function to
smooth weak discontinuities appearing in the initial data.
Accordingly, for $t<t^*$ it admits the same behavior as the standard
admissible weak solution $u$ of problem (\ref{jul767}),
(\ref{jul769}) (it can produce discontinuities in $v_{\eps}$ even
for $t<t^*$).

Therefore, for the sake of consistency (we have the same situation
for $\eps>0$ and as $\eps\to 0$), we will use the function given by
(\ref{avg2956}) in the place of $u$ appearing in (\ref{jul361}).

Weak asymptotic solution to Cauchy problem (\ref{jul361}),
(\ref{jul362}) we will solve separately in five areas of $(x,t)$
plane in which the solution is certainly smooth.
Then, we will connect solutions in those domains and
prove that the function formed in that way represents weak solution
to our problem.

In order to single out those domains we substitute (\ref{avg2956}) into (\ref{jul361}) and, formally, use Leibnitz rule for derivative of product:
\begin{gather}
v_{\eps t}+g(u_{\eps})v_{\eps x}=-(g(u_{\eps})_x v_{\eps}.
\label{sep1916}
\end{gather}

\begin{remark}
\label{okt496} We say "formally" since the function $v_{\eps}$ can
have discontinuities, at least for $\eps=0$. We repeat that this is
caused by non-smoothness of the new characteristics if we use
(\ref{avg2956}) as well as standard characteristics in $x_0=a_1$ and $x_0=a_2$.
\end{remark}

The system of characteristics corresponding to (\ref{sep1916}),
(\ref{jul362}) is:
\begin{equation}
\begin{split}
\dot{X}=&g(u_{\eps}),\ \ X(0)=x_0,\\
\dot{v}_{\eps}=&-v_{\eps}(g(u_{\eps}))_x,\ \
v_{\eps}(0)=\hat{v}(x_0).
\end{split}
\label{jul761}
\end{equation} We prove global resoluteness of this ODE system for $x_0\in [a_2-\eps A \frac{a_2+a_1}{2},a_1+\eps A \frac{a_2+a_1}{2}]$.
According to the inverse function theorem it is enough to prove that along entire temporal axis we have
$$\frac{\pa X}{\pa x_0}>0.$$Denote by $J=\frac{\pa x}{\pa x_0}$ where $x=x(x_0,t,\eps)$ is the function defined by Cauchy problem (\ref{avg615}).
We have proved in the previous theorem that $J>0$ for every $t>0$.
Recall that $$u_{\eps}(x,t)=u_0(\tilde{x}_0(x,t,\eps)),$$ where
$\tilde{x}_0$ is inverse function to the function $x$ defined trough
(\ref{avg615}). From (\ref{jul761}) we have (we write below $
g'(u_0)=g'(u_0(\tilde{x}_0(x,t,\eps)))$)
$$\frac{d}{dt}\frac{\pa X}{\pa x_0}=g'(u_0)u_0'\frac{\pa \tilde{x}_0}{\pa X}\frac{\pa X}{\pa x_0}=
g'(u_0)u_0'J^{-1}\frac{\pa X}{\pa x_0}, \ \ \frac{\pa X}{\pa
x_0}\Big{|}_{t=0}=1.$$ After integrating this differential equation
with respect to the unknown function $\frac{\pa X}{\pa x_0}$ we
obtain:
$$\frac{\pa X}{\pa x_0}={\rm exp}(\int_0^t g'(u_0)u_0'J^{-1}dt')>0, \ \ t>0,$$
which implies existence of inverse function ${x}_0={x}_0(X,t,\eps)$
along entire temporal axis, which, in turn, implies global
resoluteness of problem (\ref{jul761}).

Denote by $\varphi_i^*$, $i=1,2$, solutions of the following Cauchy problems:
\begin{align*}
\dot{X}=&g(u_{\eps}),\nonumber\\
X(0)=&a_i, \ \ i=1,2.
\end{align*}

Now, we can introduce domains in which we will separately solve Cauchy problem (\ref{jul361}), (\ref{jul362}).
\begin{remark}
In the domains to be introduced, it is equivalent to say Cauchy
problem (\ref{sep1916}), (\ref{jul362}) instead of  (\ref{jul361}),
(\ref{jul362}) since in those domains problem (\ref{jul361}),
(\ref{jul362}) has smooth solution.
\end{remark}

We set \begin{gather*}
D_1=\{(x,t)| \; x< \varphi_2\}, \ \ D_2=\{
(x,t)|\; x>\varphi_1\},\\
D_3=\{(x,t)| \; \varphi_2<x< \varphi^*_2\}, \ \ D_4=\{ (x,t)|\;
\varphi_1^*<x<\varphi_1 \},\\
D_5=\{ (x,t)| \; \varphi_2^*<x<\varphi_1^*\}
\end{gather*}

On the beginning, we prove that those domains are disjunct.
Accordingly, we inspect relations between the functions $\varphi_i$
and $\varphi_i^*$, $i=1,2$. We have to prove the following fact for
every $t\in {\bf R}^+$:
\begin{gather}
\varphi_2(t,\eps) \leq \varphi_2^*(t,\eps)<\varphi_1^*(t,\eps) \leq
\varphi_1(t,\eps), \label{sep2026}
\end{gather}

First, we prove that
\begin{gather}
\varphi_2\leq \varphi_2^*. \label{sep2036}
\end{gather} In the moment $t=0$ we have
\begin{gather}
(\varphi_2)'_t=f(U_1)(B_2-B_1)+c B_1 \text{   and   }
(\varphi_2^*)'_t=g(U_1), \label{sep2056}
\end{gather} and (see (\ref{jul366}))
\begin{align*}
g(U_1)> f'(U_1)>f'(U_1)(B_2-B_1)+cB_1,
\end{align*}since $f'(U_1)>c/2$.
Using well known theorem from ODE-s ("who goes slower does not reach
further", \cite{arn}) from (\ref{sep2056}) we see that in some
neighborhood of $t=0$ we have $\varphi_2<\varphi_2^*$. Assume now
that $t_0$ is the smallest $t>0$ such that $\varphi_2=\varphi_2^*$.
In this case we have the same situation as in the moment $t=0$, i.e.
there exists neighborhood $(t_0,t_0+\delta)$ such that
$\varphi_2<\varphi_2^*$ in $(t_0,t_0+\delta)$. Continuing like this
we see that we indeed have (\ref{sep2036}).

In the completely same manner we prove that
\begin{gather}
\varphi_1^*\leq \varphi_1. \label{sep2046}
\end{gather}


It is remained to prove that:
\begin{gather}
\varphi_2^* < \varphi_1^*. \label{sep2066}
\end{gather}This directly follows from the fact that characteristics of problem (\ref{jul761}) do not intersect.
That means that relation between two characteristics remains the
same along entire time axis. Therefore,
$$\varphi_2^*=X(a_2,t,\eps)< X(a_1,t,\eps)=\varphi_1^*,$$ since $a_2<a_1$. This proves (\ref{sep2066}).

Collecting (\ref{sep2036}), (\ref{sep2046}) and (\ref{sep2066}) we
obtain (\ref{sep2026}). Very important implication of relation
(\ref{sep2026}) and the fact that for $t>t^*$
$$\lim\limits_{\eps\to 0}\varphi_i(t,\eps)=\frac{c}{2}(t-t^*)+x^*, \ \ i=1,2$$ is the following. For $t>t^*$ we have:
\begin{gather}
\lim\limits_{\eps\to 0}\varphi_i^*(t,\eps)=\frac{c}{2}(t-t^*)+x^*, \
\ i=1,2. \label{sep2076}
\end{gather} We remind that the constants $t^*$ and $x^*$ are introduced in front of Theorem \ref{t**}.


Next, we solve problem (\ref{sep1916}), (\ref{jul362}) separately in
domains $D_i$, $i=1,...,5$.

In domains $D_1$ and $D_2$ we have $u_{\eps} \equiv const.$ and
therefore the characteristics corresponding to $v_{\eps}$ there are
straight lines. More precisely, we have:
\begin{align*}
v_{\eps}(x,t)\equiv V_1, \ \ (x,t)\in D_1,\\
v_{\eps}(x,t)\equiv V_0, \ \ (x,t)\in D_2.
\end{align*}

Another two domains are:
$$D_3=\{(x,t)| \; \varphi_2<x< \varphi^*_2\}, \ \ D_4=\{ (x,t)|\; \varphi_1^*<x<\varphi_1 \}.$$
In those domains we solve the following Cauchy problems:
\begin{align*}
v_{\eps t}&+g(u_{\eps})v_{\eps x}=-(g(u_{\eps})_x v_{\eps},\\
v_{\eps}|_{x=\varphi_1}&=V_{0} \text{   (initial data for the first Cauchy problem),}\\
v_{\eps}|_{x=\varphi_2}&=V_{1} \text{   (initial data for the second Cauchy problem)}.
\end{align*}We use standard method of characteristics.
Note that in this case characteristics emanate from the lines
$x=\varphi_i$, $i=1,2$, and not from $x$ axis as usual
The system of characteristics for ones emanating from the line
$\varphi_1$ has the form:
\begin{equation}
\begin{split}
\dot{X}=&g(u_{\eps}),\nonumber\\
\dot{v}_{\eps}=&-v_{\eps}(g(u_{\eps}))_x,\nonumber\\
X(t_0)=&\varphi_1(t_0)=x_0, \ \ v_{\eps}(t_0)=V_{0}.
\end{split}
\label{sep1926}
\end{equation} and for the characteristics emanating from the line $\varphi_2$ has the form:
\begin{equation}
\begin{split}
\dot{X}=&g(u_{\eps}),\nonumber\\
\dot{v}_{\eps}=&-v_{\eps}(g(u_{\eps}))_x,\nonumber\\
X(t_0)=&\varphi_2(t_0)=x_0, \ \ v_{\eps}(t_0)=V_{1}.
\end{split}
\label{sep1926a}
\end{equation}
Global solvability of this system can be proved in the same way as
one for the system (\ref{jul761}).

Next step is to solve second equation of (\ref{sep1926}) (or analogically of (\ref{sep1926a})). We have for problem (\ref{sep1926}):
\begin{align}
\dot{v}_{\eps}&=(-g(u_{\eps}))_x v_{\eps} \implies\nonumber\\
v_{\eps}&=V_{0}{\rm exp}(-\int_0^t(g(u_{\eps}))_xdt') \implies\nonumber\\
v_{\eps}&=V_{0}{\rm exp}(-\int_0^t(\frac{dX}{dt'})_x dt') \implies\nonumber\\
v_{\eps}&=V_{0}{\rm exp}(-\int_0^t\frac{\pa }{\pa x_0}\frac{dX}{dt}\cdot \frac{\pa \tilde{x}_0}{\pa X}dt')\implies\nonumber\\
v_{\eps}&=V_{0}{\rm exp}(-\int_0^t\frac{\frac{d}{dt}\frac{\pa X}{\pa x_0}}{\frac{\pa X}{\pa x_0}} dt')\implies\nonumber\\
v_{\eps}&=\frac{V_{0}}{\frac{\pa X}{\pa x_0}}, \label{jul766}
\end{align} and, similarly, for (\ref{sep1926a}):
\begin{gather}
v_{\eps}=\frac{V_1}{\frac{\pa X}{\pa x_0}}. \label{jul766a}
\end{gather}
Previous implies:
\begin{align*}
v_{\eps}(x,t)&=V_1\frac{\pa x_{03}}{\pa x}(x,t,\eps), \ \ (x,t)\in D_3,\\
v_{\eps}(x,t)&=V_0\frac{\pa x_{04}}{\pa x}(x,t,\eps), \ \ (x,t)\in D_4.
\end{align*} where $x_{03}=x_0(X,t,\eps)=\varphi_1(t_{01})$ and
$x_{04}=x_0(X,t,\eps)=\varphi_2(t_{02})$ (for appropriate $t_{0i}$,
$i=1,2$, depending on $(X,t)$) are inverse functions to the function
$X$ determined by (\ref{sep1926}) and (\ref{sep1926a}), respectively.

Finally, we solve problem (\ref{sep1916}), (\ref{jul362}) in the domain:
$$D_5=\{ (x,t)| \; \varphi_2^*<x<\varphi_1^*\}.$$
We apply similar procedure as in the previous case. The solution in this domain is:
$$v_{\eps}(x,t)=v_0(x_0(x,t,\eps))\frac{\pa x_{05}}{\pa x}(x,t,\eps),$$ where
$x_{05}=x_0(X,t,\eps)$ is inverse function to the function $x$
determined by (\ref{jul761}) ($x_0$ restricted on $[a_2-\eps A
\frac{a_2+a_1}{2},a_1+\eps A \frac{a_2+a_1}{2}]$).

Thus, we have proved the following theorem:
\begin{theorem}
The weak asymptotic solution to problem (\ref{sep1916}),
(\ref{jul362}) is given by the formula:
\begin{gather}v_{\eps}(x,t)=\begin{cases}
V_0, \ \ , (x,t)\in D_1,\\
V_0\frac{\pa x_{03}}{\pa x}(x,t,\eps), \ \ (x,t)\in D_3,\\
v_0(x_0(x,t,\eps))\frac{\pa x_{05}}{\pa x}(x,t,\eps), \ \ (x,t)\in D_5,\\
V_1\frac{\pa x_{04}}{\pa x}(x,t,\eps), \ \ (x,t)\in D_4,\\
V_1, \ \ (x,t)\in D_2.
\end{cases}
\label{sep2086}
\end{gather}
\end{theorem}
Let us analyze the function $v_{\eps}$ more closely. We have to
inspect its behavior on the boundaries of the domains $D_i$,
$i=1,2,3,4,5$, since in the domains we know that equation
(\ref{sep1916}) is satisfied.

On the lines $\bar{D}_1\cap \bar{D}_3$ and $\bar{D}_2\cap \bar{D}_4$
the function is continuous and therefore, in the domain
$\bar{D}_1\cup D_2\cup D_3\cup \bar{D}_4$ the function $v_{\eps}$
represents admissible weak solution to (\ref{sep1916}).

On the lines $\bar{D}_3\cap \bar{D}_5=\varphi_2^*$ and
$\bar{D}_4\cap \bar{D}_5=\varphi_1^*$ the function $v_{\eps}$ can
have discontinuities (see Remark \ref{okt496} and Example
\ref{sep2016}). Therefore, we have to check if Rankine-Hugoniot
conditions are satisfied on $\varphi_i^*$, $i=1,2$ for equation
(\ref{sep1916}) given in divergent form:
\begin{gather}
v_{\eps t}+(v_{\eps}g(u_{\eps}))_x=0.
\label{sep20106}
\end{gather}
Accordingly, we have to check:
$$(\varphi_i^*)_t=\frac{[g(u_{\eps})v_{\eps}]}{[v_{\eps}]}\Big|_{x=\varphi_i^*}.$$ Since $u_{\eps}$ is continous function from here it follows:
$$(\varphi_i^*)_t=g(u_{\eps})\Big|_{x=\varphi_i^*},$$ which is exactly the definition of the function $\varphi_i^*$, $i=1,2$.

This concludes the proof that the function $v_{\eps}$ represents weak solution to problem (\ref{sep20106}), (\ref{jul362}).

\section{Weak limit of the solution}

It remains to inspect the weak limit of the weak asymptotic solution
$(u_{\eps},v_{\eps})$ of problem (\ref{jul361}), (\ref{jul362}) for
$t>t^*$ (since for $t<t^*$ we have classical solution of the
considered problem).

We have already known from Corollary 9 that for $t\geq t^*$ we have:
\begin{gather}
w-u_{\eps}(x,t)\to\begin{cases}
U_1, \ \ x<\frac{c}{2}(t-t^*)+x^*,\\
U_0, \ \ x\geq \frac{c}{2}(t-t^*)+x^*.
\end{cases}
\label{avg1636}
\end{gather}So, we have to inspect weak limit of $v_{\eps}$. ore
precisely, in this section we will prove the following theorem:
\begin{theorem}
For every fixed $t>t^*$ the function $v_{\eps}$ given by
(\ref{sep2086}) satisfies as $\eps\to 0$
\begin{align}
\label{avg1626}
v_{\eps}&(x,t)\rightharpoonup \begin{cases}
V_1, \ \ x<\frac{c}{2}(t-t^*)+x^*,\\
V_0, \ \ x\geq \frac{c}{2}(t-t^*)+x^*
\end{cases}+\\&\qquad\Big[V_1(a_2+g(U_1)t-\frac{c}{2}(t-t^*)-x^*)+V_0(\frac{c}{2}(t-t^*)+x^*-a_1-g(U_0)t)+\nonumber
\\&\qquad\qquad\qquad\qquad\int_{a_2}^{a_1}v_0(x_0)dx_0
\Big]\delta(x-\frac{c}{2}(t-t^*)-x^*),
\nonumber
\end{align}where $\rightharpoonup$ means convergence in
the weak sense with respect to the real variable.
\end{theorem}
{\bf Proof:} To begin, note that we can write function $v_{\eps}$
from (\ref{sep2086}) in the following manner:
$$
v_{\eps}(x,t)=\hat{v}(x_0(x,t,\eps))\frac{\pa x_0}{\pa x}(x,t,\eps),
$$ where
\begin{gather}
x_0(x,t,\eps)=\begin{cases}
x-g(U_1)t, \ \ (x,t)\in \bar{D}_1,\\
x_{03}^{-1}(x,t,\eps)-g(U_1)\varphi_2^{-1}(x_{03}^{-1}(x,t,\eps)), \ \ (x,t)\in D_3 \\
\ \ \text{  (here first we go by  }x_{03}^{-1} \text{  to the line }
\varphi_2 \text{  so that  }
x_{03}^{-1}(x,t,\eps)=\varphi_2(t_0)\\
\ \ \text{  and then proceed to the line  }t=0
\text{  along the straight line  }x-g(U_1)t) ,\\
x_{05}^{-1}(x,t,\eps),\ \ (x,t)\in \bar{D}_5,\\
x_{04}^{-1}(x,t,\eps)-g(U_0)\varphi_1^{-1}(x_{04}^{-1}(x,t,\eps)), \ \ (x,t)\in D_4,\\
\ \ \text{  (here first we go by  }x_{04}^{-1} \text{  to the line }
\varphi_1 \text{  so that  }
x_{04}^{-1}(x,t,\eps)=\varphi_1(t_0) \\
\ \ \text{  and then proceed to the line  }t=0
\text{  along the straight line  }x-g(U_0)t)\\
 x-g(U_0)t, \ \ (x,t)\in \bar{D}_2,
\end{cases}
\label{sep2096}
\end{gather}and
$$\frac{\pa x_0}{\pa x}(x,t,\eps)=\begin{cases}
1, \ \ (x,t)\in \bar{D}_1,\\
\frac{\pa x_{03}}{\pa x}, \ \ (x,t)\in D_3,\\
\frac{\pa x_{05}}{\pa x}, \ \ (x,t)\in \bar{D}_5,\\
\frac{\pa x_{04}}{\pa x}, \ \ (x,t)\in D_4,\\
1, \ \ (x,t)\in \bar{D}_2.
\end{cases}
$$

We take $\eta \in C^1_0({\bf R})$ and write using (\ref{sep2086}):
\begin{align*}
&\int v_{\eps}(x,t)\eta(x)
dx=\int_{-\infty}^{\varphi_2-\eps}v_{\eps}(x,t)\eta(x) dx+
\int_{\varphi_2-\eps}^{\varphi_2^*}v_{\eps}(x,t)\eta(x) dx+\\
&\int_{\varphi_2^*}^{\varphi_1^*}v_{\eps}(x,t)\eta(x) dx+
\int_{\varphi_1^*}^{\varphi_1+\eps}v_{\eps}(x,t)\eta(x)
dx+\int_{\varphi_1+\eps}^{\infty}v_{\eps}(x,t)\eta(x) dx=\\&
\int_{-\infty}^{\varphi_2-\eps}V_1\eta(x) dx+
\int_{\varphi_2-\eps}^{\varphi_2^*} V_1\frac{\pa x_{04}}{\pa
x}\eta(x) dx+
\int_{\varphi_2^*}^{\varphi_1^*}v_0(x_0(x,t,\eps))\frac{\pa
x_{05}}{\pa x}\eta(x) dx+\\&
\int_{\varphi_1^*}^{\varphi_1+\eps}V_0\frac{\pa x_{03}}{\pa x}\eta
dx+\int_{\varphi_1+\eps}^{\infty}V_0\eta(x) dx.
\end{align*} Here, we have written $\varphi_i\pm \eps$ in order to
avoid possible $\varphi_i=\varphi_i^*$.

Then, we use the change of variables $x=X(x_0,t,\eps)$ where $X$ is
inverse function of the function $x_0=x_0(X,t,\eps)$ given by
(\ref{sep2096}). We have:
\begin{align}
\label{no5}
\int v_{\eps}(x,t)\eta(x) &
dx=\int_{-\infty}^{\varphi_2-\eps}v_{\eps}(x,t)\eta(x)
dx+\\& \int_{x_0(\varphi_2-\eps,t,\eps)}^{a_2}
V_1\eta(x(x_0,t,\eps))
dx_0+\int_{a_2}^{a_1}v_0(x_0)\eta(X(x_0,t,\eps)) dx_0+\nonumber\\&
\int_{a_1}^{x_0(\varphi_1+\eps,t,\eps)}V_0\eta(x(x_0,t,\eps))
dx_0+\int_{\varphi_1+\eps}^{\infty}V_0\eta(x) dx,
\nonumber
\end{align}and we remind that:
\begin{gather*}
x_0(\varphi_1, t, \eps)=\varphi_1+\eps-g(U_0)t, \ \ x_0(\varphi_2,
t, \eps)=\varphi_2-\eps-g(U_1)t,
\end{gather*} and for $t>t^*$ we have (see (\ref{avg25}) and (\ref{avg35})):
\begin{align*}
x_0(\varphi_1+\eps, t, \eps)&\to \frac{c}{2}(t-t^*)+x^*-g(U_0)t, \ \
\eps\to 0,
\\
x_0(\varphi_2-\eps, t, \eps)&\to \frac{c}{2}(t-t^*)+x^*-g(U_1)t, \ \
\eps\to 0.
\end{align*}

Accordingly, for $t>t^*$ after letting $\eps\to 0$ we have from
(\ref{no5}) exactly (\ref{avg1626}).

This concludes the theorem. $\Box$

It remains to give a comment concerning admissibility of the
solution. Actually, it follows from assumptions on $f'$ and $g$
quoted in Theorem \ref{th3} providing:
\begin{equation}
g(U_1)<f'(U_1), \ \ f'(U_0)<g(U_0).
\label{avg1616}
\end{equation} Admissibility of the shock wave
appearing in the solution to problem (\ref{jul767}), (\ref{jul769}) implies $$f'(U_1)<c/2<f'(U_0)$$ which together with (\ref{avg1616}) implies:
$$g(U_1)<c/2<g(U_0),$$ which proves overcompressibility of the shock and $\delta$ shock wave appearing in (\ref{avg1636}) and (\ref{avg1626}).

\begin{example}
\label{sep2016} In this example we look for the classical solution
to the following problem:
\begin{align*}
&u_t+(\frac{1}{2}u^2)_x=0\\
&v_t+(2uv)_x=0
\end{align*} for initial functions given by:
\begin{align*}
u\Big{|}_{t=0}=&\begin{cases}
1, \ \ x\leq -1\\
-x, \ \ -1\leq x\leq 1,\\
-1, \ \ 1\leq x,
\end{cases}\\
v\Big{|}_{t=0}=&\begin{cases}
1, \ \ |x|\geq 1\\
x^{2/3}, \ \ -1\leq x\leq 1.
\end{cases}
\end{align*}Using previous notation on this case, we have:
\begin{gather*}
a_1=1, \; a_2=-1, \; U=1, \; u_0^0=-1\\
b=0, \ \ K=1.
\end{gather*}
Also, for the inverse function $x_0=x_0(x,t)$ we have:
\begin{equation*}
x_0=\begin{cases}
x-2t, \ \ x\leq t-1,\\
-2-\frac{(1-t)^2}{x}, \ \ t-1\leq x\leq -(t-1)^2,\\
\frac{x}{(1-t)^2}, \ \ -(t-1)^2\leq x \leq (t-1)^2,\\
2-\frac{(1-t)^2}{x}, \ \ (t-1)^2\leq x \leq 1-t,\\
x+2t, \ \ 1-t\leq x.
\end{cases}
\end{equation*} From here we can immediately compute $\frac{\pa x_0}{\pa x}$ and then
$$
\frac{\pa x}{\pa x_0}=\frac{1}{\frac{\pa x_0}{\pa x}}.
$$

The solution of the problem is:
\begin{align}
u(x,t)&= \begin{cases}
1, \ \ x>t-1,\\
\frac{x}{t-1}, \ \ t-1\leq x<-t+1,\\
-1, \ \ -t+1 \leq x.
\end{cases}\\
v(x,t)&=\begin{cases}
1, \ \ x<t-1,\\
 \frac{(1-t)^2}{x^2},\ \ t-1\leq x <-(t-1)^2,\\
 \frac{x^{2/3}}{(1-t)^{10/3}}, \ \ -(t-1)^2\leq x < (t-1)^2,\\
 \frac{(1-t)^2}{x^2}, \ \ (t-1)^2\leq x <-t+1\\
 1, \ \ x>1-t.
\end{cases}
\label{sep916}
\end{align} Notice that the function $v$ has discontinuities on the following lines:
$$x=(t-1)^2, \ \ x=-(t-1)^2.$$
 and that on that line Rankine-Hugoniot conditions are satisfied.

 Notice, also, that if instead of $x^{2/3}$ in the initial data for
 $v$ we would put $x^{2n}$ the solution of our problem would be
 continuous for $t<t^*$.
\end{example}

\end{document}